\documentclass[11pt]{article}
\usepackage{amsmath,amsthm,amssymb}

\chardef\bslash=`\\ \hfuzz1pc

\newtheorem{defn}{Definition}[section]
\numberwithin{equation}{section}

\begin{document}

\title{On trigonometric-like  decompositions of functions with respect to the cyclic group of order n}
\author{A.K.Kwa\'sniewski*, B.K. Kwa\'sniewski** \\
\\
*Institute of Computer Science, Bia{\l}ystok University\\
PL-15-887 Bia{\l}ystok, ul.Sosnowa 64, POLAND\\
 e-mail: kwandr@uwb.edu.pl\\
** Institute of Mathematics, Bia{\l}ystok
University}
\maketitle

\begin{abstract}
The cyclic group labeled family of $\alpha-$projection operators
implicitly present in  [23] is  used  as in [4-8,16] for
investigation of  decomposition of functions with  respect to the
cyclic group of order n.    Series of new identities thus arising
are demonstrated and new perspectives for further investigation
are indicated as for example in the case of q-extended special
polynomials.   The paper constitutes an example of the application
of the method  of projections introduced  in  [21] ; see also
references [4-8] .

\end{abstract}

{\small {\ KEY WORDS: cyclic group labeled projection operators,
special hyperbolic functions }

\section{\protect Introduction}
{\small In the past century Ungar had introduced in  his  Indian
J. Pure Appl. Math paper [23]  higher order $\alpha-$hyperbolic
functions which  are denoted here as $\ Z_{n}$ cyclic group
labelled family $\{ \ h_{s}^{\alpha}(z)\}_{s\in Z_{n}}$. These
functions  are specific examples of  eigenfunctions of the scaling
$\Omega$ operator. $\Omega$  operator is used in this note to
define a family of mutually  orthogonal $\alpha-$projection
operators $\{\Pi_{i}^{\alpha}\}_{l\in Z_{n}}$ and then these
eigenfunctions themselves - including Ungar`s $\alpha-$hyperbolic
functions (see: further examples below). For more information on
the history of rediscovering standard  $\alpha=1-$hyperbolic
$\{h_{s}(z)\}_{s\in Z_{n}}$ and $\alpha=-1-$hyperbolic i.e.
circular  $\{f_{s}(z)\}_{s\in Z_{n}}$ functions see first of all
\textit{The Mathematics Magazine} article [21] and for further
references see also [14] and [22]. If one takes $n=2$, and
$\alpha=\pm1$ then  one obtains  cosh , sinh or cos and sin
functions. Here one proposes an expedition a little bit more
farther then \textit{"Beyond Sin and Cos ?"}[21]. In our story
$Z_{n}=\{0,1,...,n-1\}$ denotes cyclic group under the addition
i.e. for $k,l\in Z_{n}$: $k\buildrel {\cdot} \over +l$ denotes
addition mod n and $k\buildrel {\cdot}\over - l$ denotes
subtraction mod n; $\omega=\exp(i\frac{2\pi}{n})$; $n>1$. While
extending the range of \textit{"Beyond Sin and Cos"}[21] we use
group $Z_{n}$ labeled family of $\alpha$-projection
$\{\Pi_{l}^{\alpha}\}_{l\in Z_{n}}$ operators. Projection
operators $\{\Pi_{l}\}_{l\in Z_{n}}$ i.e.
$\{\Pi_{l}^{\alpha}\}_{l\in Z_{n}}$ with $\alpha=1$ were used in
[4-8,16] for investigation of decomposition of functions with
respect to the cyclic group of order n. We arrive at the
\textit{"Beyond Sin and Cos"} while $\alpha$-decomposing $\exp$
function. Does then decomposing from [4,16] of a function L given
by Laurent series L leads too far \textit{"Far Beyond Sin and
Cos?"} Perhaps this would be the better title of this article.}

\section{\protect Expected Elementary Background}
{\small Apart from \textit{The Mathematics Magazine} article [21]
also monograph [9] on circulant matrices is recommended. As for
Laurent series L considered here these may be treated also as
formal Laurent series. This includes algebras of formal series
(formal power series, exponential formal power series, Dirichlet
series  etc.) - as used in combinatorics [24]. This aspect is not
pursued here - let us however remark that projection operators
$\{\Pi_{l}\}_{l\in Z_{n}}$ are ready to be applied for a might be
desirable study of $Z_{n}$ labeled subsequences of counting
sequences in combinatorics. The use of circulant matrices enables
one to introduce $Z_{n}$- L- correspondents of trigonometric
formulas for hyperbolic $\{h_{s}(z)\}_{s\in Z_{2}}$ functions of
second order in a manner this was done for hyperbolic
$\{h_{s}(z)\}_{s\in Z_{n}}$ functions of $n$-th order in [22,14]}.

{\small It is easy to see that for $\{h_{s}(z)\}_{s\in
Z_{2}}\equiv\{\cosh z,\sinh z\}\ z\in \mathcal{C}$ from the  group
property $$\left(\begin{array}{cc}
  \cosh z & \sinh z\\
  \sinh z & \cosh z
\end{array}\right)
\left(\begin{array}{cc}
  \cosh w & \sinh w\\
  \sinh w & \cosh w
\end{array}\right)=
\left(\begin{array}{cc}
  \cosh(z+w) & \sinh(z+w)\\
  \sinh(z+w) & \cosh(z+w)
\end{array}\right) \forall w,z \in{\mathcal{C}}$$

de Moivre formulas in their matrix form follow:
$$\left(\begin{array}{cc}
  \cosh z & \sinh z\\
  \sinh z & \cosh z
\end{array}\right)^{n}=\left(\begin{array}{cc}
  \cosh nz & \sinh nz\\
  \sinh nz & \cosh nz
\end{array}\right) \forall w,z \in{\mathcal{C}}$$.
\\$Z-{n}$ - exp-counterparts of the hyperbolic-trigonometric
identity: $\forall w,z \in{\mathcal{C}}$; $$ (\cosh z)^{2}-(\sinh
z)^{2}=1 i.e. \det \left(\begin{array}{cc}
  \cosh z & \sinh z\\
  \sinh z & \cosh z
\end{array}\right)=1$$
\\is introduced with help of  circulants as in [14] (see also [2]).
In the sequel we shall try to answer the question: are there
$Z_{n}-$ L-counterparts also available?

The set of matrices $\left(\begin{array}{cc}
  \cosh z & \sinh z\\
  \sinh z & \cosh z
\end{array}\right)=1$; $z\in R$ under matrix multiplication
constitutes SO(1.1) group. This is the group of two dimensional
special relativity transformations. The set of matrices
$\left(\begin{array}{cc}
  \cos z & -\sin z\\
  \sin z & \cos z
\end{array}\right)$;  $z\in R$ under matrix multiplication  constitutes SO(2) group;  this is of
course the group of  two dimensional rotations.}

\section {\protect $Z_{n}$ cyclic group labeled $\alpha$-projection
operators and $Z_{n}$ decomposition of functions}

{\small In this section $Z_{n}$ labeled  $\alpha$-projection
operators are used for decomposition of functions with  respect to
the cyclic group of order $n$ [23,21], [4-8,16].  Let us then
define this family of  $\alpha$-projection operators.}

\begin{defn}
{\small $\{\Pi_{l}^{\alpha}\}_{l\in Z_{n}}$ acting on the linear
space of functions of complex variable are defined  according to
\begin{equation}\label{e1}
\Pi_{k}^{\alpha}:=\frac{1}{n}\alpha^{-\frac{k}{n}}\sum\limits_{s\in
Z_{n}}{\omega^{-ks}\Omega^{s}S(\sqrt[n]{\alpha})}
\end{equation}
where $\sqrt[n]{\alpha}$ is an arbitrarily specified n-th root of
$\alpha$ and $\Omega$, $S(\lambda)$ are scaling operators:
\begin{equation}\label{e2}
(\Omega f)(z):=f(\omega z), \mbox{  }(S(\lambda )f)(z):=f(\lambda
z) \mbox{  }i.e. \Omega = S(\omega).
\end{equation}
}\end{defn}
The family of $\alpha $-projection operators
$\{\Pi_{l}^{(\alpha)}\}_{l\in Z_{n}}$ extends the set of families
of projection operators $\{V_{k}\}_k\in Z_{n}$; $V_{k} \cdot
V_{l}=V_{l}\delta_{kl}$ introduced in [15].
$\{\Pi_{l}^{(\alpha)}\}_{l\in Z_{n}}$ is an easy generalization of
the family of projection operators used under notation
$\{\Pi_{[n,k]}\}_{k\in Z_{n}}$ in  [4-8] for a decomposition of
various special functions with respect to the cyclic group of
order n  in analogy to the  decomposition of $\exp$ function
standard hyperbolic functions of  n-th  as was done and used under
the notation $\{\Delta_{k}\}_{k \in {Z_{n}}}$ in [16] in order to
investigate higher order recurrences for analytical functions of
Tchebysheff type [16,1]. As $\{\Delta_{k}\}_{k \in {Z_{n}}}\equiv
\{\Pi_{[n,k]}\}_{k\in Z_{n}}$ we shall use notation
$\{\Pi_{k}\}_{k\in Z_{n}}\equiv \{\Pi_{[n,k]}\}_{k\in Z_{n}}
\equiv \{\Pi_{k}^{\alpha =1}\}_{k\in Z_{n}}\equiv
\{\Delta_{k}\}_{k \in {Z_{n}}}$ in conformity with all the papers
mentioned and also this note. Of course (one arguments like in
[15]) $\Pi_{l}\Pi_{m}=\delta_{lm}\Pi_{l}$ and from (\ref{e1}) and
(\ref{e2}) one sees that
$\Pi_{k}^{(\alpha)}=\alpha^{\frac{-k}{n}}S(\sqrt[n]{\alpha})\Pi_{k}$.
Hence we infer what follows.\\
\textbf{Observation 3.1}
\begin{equation}\label{e3}
\Pi_{l}^{(\alpha)}\Pi_{m}^{(\alpha)}=\delta_{lm}\Pi_{l}^{(\alpha)}\alpha^{\frac{-m}{n}}S(\sqrt[n]{\alpha})
\end{equation}
\begin{equation}\label{e4}
\sum\limits_{k\in{Z_{n}}}\alpha^{\frac{k}{n}}\Pi_{k}^{(\alpha)}=S(\sqrt[n]{\alpha})
\mbox{  } end \sum_{k\in {Z_{n}}} \Pi_{k}=id.
\end{equation}
Although - as seen from formulas (\ref{e1}),(\ref{e2}),(\ref{e3})
and (\ref{e4})-the $\alpha$-projection operators
$\Pi_{k}^{(\alpha)}$ differ from projection operators $\Pi_{k}$
only by rescaling we keep introducing them because of reasons
$\alpha$-hyperbolic functions were introduced in [21,23]. Namely
$\alpha=-1, 0, +1$ cases may be treated with the same method and
then formulas specified. This will therefore include $\alpha=-1-$
hyperbolic i.e. circular $\{f_{s}(z)\}_{s\in Z{n}}$ functions,
$\alpha=0$ - hyperbolic i.e. "binomial" [21]
$\{h_{s}^{0}(z)=\frac{z^{s}}{s!}\}_{s\in Z_n}$ functions and
$\alpha=1-$ hyperbolic i.e. hyperbolic $\{h_{s}(z)\}_{l\in Z_{n}}$
functions. (In the $\alpha=0$ case one uses after [23,21] the
convention $0^{0}=1$- see Example 3.1.). Moreover, with help of
these $\alpha$-projection operators $\{\Pi_{l}^{(\alpha)}\}_{l\in
Z_{n}}$ one may define new families of eigenfunctions of the
$\Omega$ operator. Here there are some introductory examples based
on [16]. \\
\textbf{Example 3.1.} \\
Let $\{h_{s}^{(\alpha)}(z)\}_{s\in {Z_{N}}}$ where
$h_{s}^{\alpha}:=\Pi_{s}^{(\alpha)}\exp$ then
$h_{s}^{\alpha}=\sum\limits_{k\geq
0}\frac{\alpha^{k}z^{nk+s}}{(nk+s)!}=\alpha^{\frac{-s}{n}}\sum\limits_{k\geq
0}{\frac{(\sqrt[n]{\alpha}z)^{nk+s}}{(nk+s)!}}$ and $\Omega
h_{s}^{\alpha}=\omega^{s}h_{s}^{\alpha}$; $s\in {Z_{n}}$. We shall
call: $h_{l}^{\alpha}$ the $l-\alpha$-hyperbolic series (compare
with [23,21]). Of course
\begin{equation}\label{e5}
h_{s}^{\alpha}=\frac{1}{n}\alpha^{\frac{-s}{n}}\sum\limits_{k\in{Z_{n}}}\omega^{-ks}\exp(\omega^{k}\sqrt[n]{\alpha}z)
\end{equation}

Note also [14] for future use that for $h_{l}\equiv
h_{l}^{\alpha=1}$
\begin {equation}\label{e6}
\exp(\omega_{l}z)=\sum\limits_{k\in{Z_{n}}}\omega^{kl}h_{k}(z)
\end{equation}
Let $\{g_{l}^{\alpha}(z)\}_{l\in{Z_{n}}}$ where
$g_{l}^{\alpha}:=\Pi_{l}^{(\alpha)}\frac{1}{1-id}$ with
$\frac{1}{1-id}(z)\equiv \frac{1}{1-z}$ and $l\in Z_{n}$ then
$g_{l}^{\alpha}(z)=\sum\limits_{k\geq 0}\alpha^{k}z^{nk+1}$ and
$\Omega g_{l}^{\alpha}=\omega ^{l}g_{l}^{\alpha}$. We shall call:
$g_{l}$ the $l-\alpha$-geometric series; (compare with [16]). Of
course \begin{equation}\label{e7}
g_{s}^{\alpha}(z)=\frac{1}{n}\alpha^{\frac{-s}{n}}\sum\limits_{k\in{Z_{n}}}\omega^{-ks}g(\omega^{k}\sqrt[n]{\alpha}z)
\end{equation}
Let $\{L_{l}^{\alpha}(z)\}_{l\in Z_{n}}$ where
$L_{l}^{\alpha}:=\Pi_{l}^{(\alpha)}L$ with $L(z)=\sum\limits_{k\in
Z}a_{k}z^{k}$ and $l\in Z_{n}$ then
$L_{l}^{\alpha}=\sum\limits_{k\in Z}a_{nk+l}\alpha^{k}z^{nk+l}$
and $\Omega L_{l}^{\alpha}=\omega^{l}L_{l}^{\alpha}$. We shall
call: $L_{l}^{\alpha}$ the $l-\alpha$-Laurent series; (compare
with [16]). Of course
\begin{equation} \label{e8}
L_{s}^{\alpha}(z)=\frac{1}{n}\alpha^{\frac{-s}{n}}\sum\limits_{k\in
Z_{n}}\omega^{-ks}L(\omega^{k}\sqrt[n]{\alpha}z)
\end{equation}
Indeed: $\Pi_{l}\sum\limits_{k\in Z}a_{k}z^{k}=\sum\limits_{k\in
Z}a_{k}\frac{1}{n}\sum\limits_{s\in
Z_{n}}\omega^{s(k\dot{-}l)}z^{k}=\sum\limits_{k\in
Z}a_{k}z^{k}\delta(k\dot{-}l)=\sum\limits_{m\in
Z}a_{nm+l}z^{mn+l}$ and now act with $\alpha
\frac{-1}{n}S(\sqrt[n]{\alpha})$ on both sides in order to see
that $L_{l}^{\alpha}(z)=\sum\limits_{k\in
Z}a_{nk+l}\alpha^{k}z^{nk+l}$. This simple method of decomposition
([23,21] and [4-8,16]) of functions with respect to $Z_{n}$ just
by acting on them as in this note by $\alpha$ projection operators
$\{\Pi_{l}^{(\alpha)}\}_{\in Z_{n}}$ may be then extended to
explore special properties of new special functions
$L_{l}^{\alpha}(z)$; $L\in \textbf{Z}_{n}$ where $L$ is any
function expandable around complex $0\in \mathcal{C}$ into Laurent
series. In view of {(\ref {e4})} functions $L_{l}(z)\equiv
L_{l}^{\alpha =1}(z)$; $l\in \textbf{Z}_{n}$ "\textit{preserve the
flavour of striking results like Euler's formula}" [21]. Indeed -
in our $\alpha=1$ case the generalized Euler formula is just this:
\begin{equation} \label{e9}
\sum\limits_{k\in Z_{n}}L_{l}(z)=L(z).
\end{equation}
Also analogue of de Miovre formulas presented here in their matrix
form [14] holds as well as correspondents of $(\cosh z)^{2}-(\sinh
z)^{2}=1$. \\
\textbf{Example 3.2} $n=3; \alpha=1$ case. \\Let us consider
generalizations of cosh and sinh hyperbolic functions of the
second order known since a long time (see [21,22],[14]). They are
defined according to ($h_{l}\equiv h_{l}^{\alpha=1}$)
\begin{equation} \label{e10}
h_{i}(x)=\frac{1}{3}\sum\limits_{k\in
Z_{3}}\omega^{-ki}\exp\{\omega^{k}x\}; i\in \textbf{Z}_{3};
\omega=\exp\left\{i\frac{2\pi}{3}\right\}
\end{equation}
One may call also (\ref{e10}) - Euler's formulas for hyperbolic
functions of $n$-th order with $n=3$. We put $n=3$ only for
convenience of easy presentation. In [14] identities for
$\{h_{i}\}_{i\in \textbf{Z}_{n}}$ hyperbolic functions are derived
from properties of "de Moivre" groups which for $z\in R$ and $m=2$
coincide with $SO(1,1)$ (hyperbolic case). Let us then introduce
at first $\gamma=(\delta_{i,k\dot{-1}})$; $k,i\in \textbf{Z}_{3}$
i.e. $\gamma= \left(\begin{array}{ccc}
  0 & 1 & 0 \\
  0 & 0 & 1 \\
  1 & 0 & 0
\end{array}\right)$ - the matrix generator of this "de Moivre"
one parameter group. Now the following is obvious (check it):
$\gamma^{n}=(\delta_{i,k\dot{-}i})^{n}=I$ and
$Tr\gamma=Tr(\delta_{i,k\dot{-}i})=0$. Hence $\det\exp\{\gamma
z\}=\exp\{Tr\gamma z\}=1$ and $\{H(z)=\exp\{\gamma z\}\}_{z\in C}$
forms what we call de Moivre group because $H(z)H(w)=H(z+w)$ where
$H(z)=\exp\{\gamma z\}$; $\gamma=(\delta_{i,k\dot{-}i})$; $k,i\in
\textbf{Z}_{n}$ and $\det H(z)=1$. Thus we arrive at the following
observation. \\\textbf{Observation 3.2} \textit{de Moivre formulas
for $n=3$ in their matrix may by written as follows: $\forall \phi
\in C$ and $\forall n\in Z$}
\begin{equation}\label{e11}
  H(n\phi)=\left (\begin{array}{ccc}
    h_{0}(n\phi) & h_{1}(n\phi) & h_{2}(n\phi) \\
    h_{2}(n\phi) & h_{0}(n\phi) & h_{1}(n\phi) \\
    h_{1}(n\phi) & h_{2}(n\phi) & h_{0}(n\phi) \
  \end{array}\right)
\end{equation}
Due to the group property of $\{H(z)=\exp\{\gamma z\}\}_{z\in C}$
one easily gets series of identities [14,16].\\\textbf{Observation
3.3} \textit{For $n=3$: $\forall k,m\in Z$ and $\forall l\in
Z_{3}$ the following three identities hold:}
\begin{equation}\label{e12}
  h_{i}((n+k)z)=3h_{0}(nz)h_{l}(kz)-h_{l}((n+k\omega)z)-h_{l}((n+k\omega^{2})z)
\end{equation}
\textit{as well as $\forall x\in C$}
\begin{equation}\label{e13}
  h_{0}(3x)=h_{0}^{3}(x)+h_{1}^{3}(x)+h_{2}^{3}(x)+3!h_{0}(x)h_{1}(x)h_{2}(x)
\end{equation}
and (see[17])
\begin{equation}\label{e14}
  h_{0}(x)h_{1}(x)h_{1}(x)=\frac{1}{9}(h_{0}(3x)-1)
\end{equation}
\\\textbf{Observation 3.4} \textit{The identity corresponding to
$(\cosh\alpha)^{2}-(\sinh\alpha)^{2}=1$ identity for $n=2$ is the
following:}
\begin{equation}\label{e15}
  h_{0}^{3}(\phi)+h_{1}^{3}(\phi)+h_{2}^{3}(\phi)-3h_{0}(\phi)h_{1}(\phi)h_{2}(\phi)=1
\end{equation}
\textit{which is equivalent to $\det H(\phi)=1; \forall\phi\in C$
i.e.}
\begin{equation}\label{e16}
  \left|\begin{array}{ccc}
    h_{0}(\phi) & h_{1}(\phi) & h_{2}(\phi) \\
    h_{2}(\phi) & h_{0}(\phi) & h_{1}(\phi) \\
    h_{1}(\phi) & h_{2}(\phi) & h_{0}(\phi) \
  \end{array}\right|=1
\end{equation}
Naturally some of these identities are easy to be written for
arbitrary $n$; for example  (\ref{e12}) is a specification of
(\ref{e17}) (see: (2.1) in [16]).\\ \textbf{Observation 3.5.}
$\forall \alpha,\beta\in C$ and $\forall l\in Z_{n}$:
\begin{equation}\label{e17}
  h_{0}(\alpha)h_{l}(\beta)=\frac{1}{n}\sum\limits_{k\in Z_{n}}h_{l}(\alpha
  + \omega ^{k}\beta)
\end{equation}
(Compare the formula (\ref {e17}) with (4.2) formulas in [22].)\\
From $H(z)H(w)=H(z+w)$ and cyclicity of $H(z)$ matrix we derive:
$\forall k\in Z_{n}$
\begin{equation}\label{e18}
  h_{k}(x+y)=\sum\limits_{i\in Z_{n}}h_{i}(x)h_{k\dot{-}i}(y)
\end{equation}
\textbf{For real} parameter $\phi\in R$ the elements of the de
Moivre one parameter group might be represented  by points
$(h_{0}(\phi),h_{1}(\phi),h_{2}(\phi))$ of the curve defined by
(\ref {e16}). This curve runs on the surface defined by the
equation $x^{3}+y^{3}+z^{3}-3xyz=1$; see [3]. For $m=4$ case due
to $\det H(\phi)=1$ - the corresponding hyper-surface is defined
by
$$-x^{4}+y^{4}-z^{4}+t^{4}+4x^{2}yt-4xy^{2}z+4z^{2}yt-4t^{2}xz+2x^{2}z^{2}-2y^{2}t^{2}=1$$
\textbf{Example $\alpha-$3.2.} $n=3$\textit{ - Ungar`s} $\alpha-$
\textit{hyperbolic case}. \\Let us consider Ungar`s $\alpha$-
hyperbolic functions of the $n=3$ order. They are defined - by
(\ref {e5})
$$h_{i}^{\alpha}(z)=\frac{1}{3}\alpha^{\frac{-i}{n}}\sum\limits_{k\in
\textbf{Z}_{3}}\omega^{-ki}\exp\{\omega^{k}\sqrt[n]{\alpha}z\};
i\in \textbf{Z}_{3}; \omega=\exp\{i\frac{2\pi}{3}\}$$ One may call
$(\alpha-10)$ - Euler's formulas for Ungar`s  $\alpha-$hyperbolic
functions of $n-$th order with $n=3$.(We put $n=3$ only for
convenience of easy presentation). As in [16] identities for
$\{h_{s}^{\alpha}(z)\}_{s\in \textbf{Z}_{n}}$ Ungar`s
$\alpha-$hyperbolic functions might be derived from properties of
$\alpha-$de Moivre  groups which for $z\in R$ and $m=2$ coincide
with $SO(2)$ ($\alpha=-1$; elliptic case) or $SO(1,1)$ ($\alpha=
+1$; hyperbolic case). The matrix generator of this $\alpha-$de
Moivre one parameter group is the matrix
$$\gamma(\alpha)=\left(\begin{array}{ccc}
  0 & 1 & 0 \\
  0 & 0 & 1 \\
  \alpha & 0 & 0
\end{array}\right).$$ The following is obvious: $\gamma(\alpha)^{n}=\alpha I$, $Tr\gamma=0$
and as $\det\{\exp A\}=\exp\{TrA\}$ then
$\det\exp\{\gamma(\alpha)z\}=1$ and
$\{H^{\alpha}(z)=\exp\{\gamma(\alpha)z\}\}_{z\in C}$ forms what we
call an $\alpha-$de Moivre group. Sure;
$H^{\alpha}(z)H^{\alpha}(w)=H^{\alpha}(z+w)$ takes place for
arbitrary $n$ where $H^{\alpha}(\phi)=\exp\{\gamma(\alpha)\phi\}$,
$\gamma(\alpha)=(\delta_{i,k\dot{-i}}+(\alpha-1)\delta_{n-1,0})$;
$k,i\in \textbf{Z}_{n}$ and $\det H^{\alpha}(z)=1$. Therefore we
observe what follows. \\\textbf{Observation 3.6.}
$\alpha-$\textit{de Moivre formulas in the matrix form are given
by:}
$$\forall \alpha,\phi\in C, \forall n\in Z, H^{\alpha}(n\phi)=\left(\begin{array}{ccc}
  h_{0}^{\alpha}(n\phi) & h_{1}^{\alpha}(n\phi) & h_{2}^{\alpha}(n\phi) \\
  \alpha h_{2}^{\alpha}(n\phi) & h_{0}^{\alpha}(n\phi) & h_{1}^{\alpha}(n\phi) \\
  \alpha h_{1}^{\alpha}(n\phi) & \alpha h_{2}^{\alpha}(n\phi) & h_{0}^{\alpha}(n\phi)
\end{array}\right)\equiv (H^{\alpha}(\phi))^{n}$$
\textbf{For real} group parameter $\phi \in R$ and $\alpha\in R$
the elements of the de Moivre one parameter group might be
represented by points
$(h_{0}^{\alpha}(\phi),h_{1}^{\alpha}(\phi),h_{2}^{\alpha}(\phi))$
of the curve defined by $\det H^{\alpha}(\phi)=1$. This curve runs
on the surface defined by the equation
$$x^{3}+\alpha y^{3}+\alpha^{2}z^{3}-\alpha 3xyz=1.$$
Due to the group property of
$\{H^{\alpha}(\phi)=\exp\{\gamma(\alpha)\phi\}\}_{z\in C}$ one may
obtain series of identities as in [14,16]. \\\textbf{Problem 3.1.}
\textit{We ask: "May one obtain trigonometric-like identities
(\ref {e10})-(\ref{e18}) in the case when exp function is replaced
by L function representing Laurent series"?}\\For that to try to
answer in $\alpha=1$ case let us at first recall again the
identity (\ref {e6}) and let us note that it generalizes to the
case when exp function is replaced by L function representing
Laurent series-just act with $\sum\limits_{k\in Z_{n}}\Pi_{k}=id$
on $L(\omega^{l}z)$ - i.e.
\begin{equation}\label{e19}
L(\omega^{l}z)=\sum\limits_{k\in
Z_{n}}L_{k}(\omega^{l}z)=\sum\limits_{k\in
Z_{n}}\omega^{kl}L_{k}(x).
\end{equation}
Now introduce  circulant matrix - an analogue of $H(z)$ -
according to
\begin{equation}\label{e20}
C(\vec L)(z)=L\{\gamma z\}
\end{equation}
where $\vec L\equiv (L_{0}(z),L_{1}(z),\dots,L_{n-1}(z))$ i.e.
consider the matrix of the form

\begin{equation}\label{e21}
C(\vec L)(z)=\left(\begin{array}{cccc}
  L_{0}(z) & L_{1}(z) & \dots & L_{n-1}(z) \\
  L_{n-1}(z) & L_{0}(z) & \dots & L_{n-2}(z) \\
  \dots & \dots & \dots & \dots \\
  L_{1}(z) & L_{2}(z) & \dots & L_{0}(z)
\end{array}\right).
\end{equation}
We know from textbooks [18] that $C(\vec L)(z)=\sum\limits_{k\in
Z_{n}}L_{k}(z)\gamma^{k}$ and the spectrum of $\gamma$ matrix is
just the multiplicative cyclic group $\hat
Z_{n}=\{\omega^{k}\}_{k\in Z_{n}}$. Due to this elementary fact
\begin{equation}\label{e22}
  \det C(\vec L)(z)=\prod\limits_{l\in Z_{n}}\sum\limits_{k\in
  Z}L_{k}(z)\omega^{kl}.
\end{equation}
(see Remarks 3.1. bellow for links with discrete Fourier
transform). Formulas (\ref{e19}) and  (\ref{e22}) imply then
another one and very important one (compare with (8.1) in [5]):
\begin{equation}\label{e23}
det\left(\begin{array}{cccc}
  L_{0}(z) & L_{1}(z) & \dots & L_{n-1}(z) \\
  L_{n-1}(z) & L_{0}(z) & \dots & L_{n-2}(z) \\
  \dots & \dots & \dots & \dots \\
  L_{1}(z) & L_{2}(z) & \dots & L_{0}(z)
\end{array}\right)=\prod\limits_{l\in Z_{n}}L(\omega^{l}z)
\end{equation}

For $L=\exp$  we come back to  $\{L_{s}(z)\}_{s\in
Z_{n}}=\{h_{s}(z)\}_{s\in Z_{n}}$ i.e. hyperbolic functions of
$n-$th order and the formula (\ref {e23}) coincides with $\det
H(z)=1$ because $\prod\limits_{l\in Z_{n}}\exp(\omega^{l}z)\equiv
1$.\\
\textbf{Answer to the Problem 3.1. }: Coming now back to our
question above for $\alpha=1$ case "\textit{May one obtain
trigonometric-like identities (\ref{e10})-(\ref{e18}) in the case
when $exp$ function is replaced by L function representing Laurent
series}"? we answer: \\\textbf{Observation 3.7.} \textit{The
crucial trigonometric-like identities
(\ref{e11}),(\ref{e16})-(\ref{e18}) do not hold in the case when
exp function is replaced by  $L\neq\exp$ function representing
Laurent series"} We readily see  that  exp  function  is
exceptional and irreplaceable  because of the following.
\\\textbf{Observation 3.8.}  \textit{Only for $L=\exp$ (up to scaling of the
argument) circulant matrices $C(\hat L)(z)=L\{\gamma z\}$ form a
group such that $L(z)L(w)=L(z+w)$}.
\\However not all is lost. One may find out many counterparts,
analogue identities  to those originating from $\exp$
decomposition with help of projection operators family
$\{\Pi_{i}\}_{l\in Z_{n}}$ [23, 21,4-8,16] even in arbitrary
$\alpha\in C$ case. For that to see let us recall again  the
formula (\ref{e19}):
$$L(\omega^{l}z)=\sum\limits_{k\in
Z_{n}}L_{k}(\omega^{l}z)=\sum\limits_{k\in
Z_{n}}\omega^{kl}L_{k}(x).$$ Now apply to both sides the operator
$\sum\limits_{k\in
Z_{n}}\alpha^{\frac{k}{n}}\Pi_{k}^{(\alpha)}=S(\sqrt[n]{\alpha})$
and recall that $L_{l}^{\alpha}:=\Pi_{l}^{(\alpha)}L$. Then we get
\begin{equation}\label{e19a}
  L(\omega^{l}\sqrt[n]{\alpha}z)=\sum\limits_{k\in
  Z_{n}}\alpha^{\frac{k}{n}}L_{k}^{\alpha}(\omega^{l}z)=\sum\limits_{k\in
  Z_{n}}\alpha^{\frac{k}{n}}\omega^{kl}L_{k}^{\alpha}(x)
\end{equation}

One may also introduce $\alpha-$circulant matrix - an analogue of
$H^{\alpha}(z)$ - according to
\begin{equation}\label{e20a}
C^{\alpha}(\vec L)(z)=L\{\gamma(\alpha)z\}
\end{equation}
where $\vec L\equiv
(L_{0}^{\alpha}(z),L_{1}^{\alpha}(z),\dots,L_{n-1}^{\alpha}(z))$
so that we consider now the $\alpha-$circulant matrix

\begin{equation}\label{e21a}
C^{\alpha}(\vec L)(z)=\left (\begin{array}{cccc}
  L_{0}^{\alpha}(z) & L_{1}^{\alpha}(z) & \dots & L_{n-1}^{\alpha}(z) \\
  \alpha L_{n-1}^{\alpha}(z) & L_{0}^{\alpha}(z) & \dots & L_{n-2}^{\alpha}(z) \\
  \dots & \dots & \dots & \dots \\
  \alpha L_{1}^{\alpha}(z) & \alpha L_{2}^{\alpha}(z) & \dots & L_{0}^{\alpha}(z)
\end{array} \right)
\end{equation}
We know from textbooks [18] that $C^{\alpha}(\vec
L)(z)=\sum\limits_{k\in Z_{n}}L_{k}^{\alpha}(z)\gamma(\alpha)^{k}$
and the spectrum of $\gamma(\alpha)$ matrix is just
$\alpha^{\frac{1}{n}}\hat Z_{n}\equiv
\left\{\alpha^{\frac{1}{n}}\omega^{k}\right\}_{k\in Z_{n}}$
because $\gamma(\alpha)^{n} =\alpha I$. Due to this simple fact
and (\ref{e19a})
\begin{equation}\label{e22a}
\det C^{\alpha}(\vec L)(z)=\prod\limits_{l\in
Z_{k}}\sum\limits_{k\in
Z_{n}}L_{k}^{\alpha}(z)\alpha^{\frac{k}{n}}\omega^{kl}.
\end{equation}
\begin{equation}\label{e23a}
 \det \left (\begin{array}{cccc}
  L_{0}^{\alpha}(z) & L_{1}^{\alpha}(z) & \dots & L_{n-1}^{\alpha}(z) \\
  \alpha L_{n-1}^{\alpha}(z) & L_{0}^{\alpha}(z) & \dots & L_{n-2}^{\alpha}(z) \\
  \dots & \dots & \dots & \dots \\
  \alpha L_{1}^{\alpha}(z) & \alpha L_{2}^{\alpha}(z) & \dots & L_{0}^{\alpha}(z)
\end{array} \right)=\prod\limits_{l\in
Z_{n}}L(\omega^{l}\sqrt[n]{\alpha}z)
\end{equation}
- so as we see -  (\ref {e23a})  for $\alpha \neq 1$  is also
comfortable and handy as (\ref {e23}) (see (8.1) in [5]). For
$L=\exp$ we come back to $\{L_{s}(z)\}_{s\in
Z_{n}}=\{h_{s}^{\alpha}(z)\}_{s\in Z_{n}}$ Ungar`s -hyperbolic
functions of $n-$th order and the formula (\ref {e23a})coincides
with $\det H^{\alpha}=1$ because $\prod\limits_{l\in
Z{n}}\exp(\omega^{l}|\alpha|z)\equiv 1$.
\\\textbf{Miscellaneous Remarks 3.1.}
\begin{enumerate}
  \item $\gamma$ matrix plays a crucial role in  $Z_{n}-$quantum mechanics (see: (2.3) in [19] and
        references therein and also see: (2.5) in [19]).
  \item Columns of Sylvester matrix $S=\frac{1}{\sqrt{n}}(\omega^{kl})_{kl\in Z_{n}}$ are eigenvectors of $\gamma$ matrix
which makes a link to $\textbf{Z}_{n}-$group discrete Fourier
transform analysis and synthesis [17], where harmonic analysis is
the passage from functional values to coefficients while harmonic
synthesis is the passage from coefficients to functional values.
  \item For a related simple generalization of analytic function
theory see [11]
\item It is obvious that $H^{\alpha}(z)=\exp\{\gamma(\alpha)z\}$,
 $\gamma(\alpha)=(\delta_{i,k\dot{-}1}+(\alpha-1)\delta_{n-1,0})$; $k,i\in\textbf{Z}_{n}$
  is the unique solution of
the equation
$\frac{d}{dz}H^{\alpha}(z)=\gamma(\alpha)H^{\alpha}(z)$with
$H^{\alpha}(0)=I$. This is equivalent to say that $\gamma$ is the
generator of $\alpha-$de Moivre group
$H^{\alpha}(z)=\exp\{\gamma(\alpha)z\}$. Naturally from the above
we conclude (compare with (5) in [21]) that
$$\frac{d^{n}}{dz^{n}}H^{\alpha}(z)=\alpha H^{\alpha}(z);
\frac{d}{dz}h_{s}^{\alpha}(z)=(1+(\alpha-1)\delta_{0,s})
h_{s-1}^{\alpha}(z); s\in Z_{n}.$$
  \item Sylvester matrix of  $\textbf{Z}_{n}-$ discrete Fourier
transform analysis serves to diagonalize our hero-circulant matrix
 $C(\vec L)(z)$ as well as $\alpha-$hero : an $\alpha-$circulant matrix (see
definition  below and [21,23]).
\item Consider $\alpha\in C$ case. It is like we went away
too far from the source of trigonometric analogies i.e. from  exp
function when  exp  function is replaced by $L$ function
representing Laurent series. Therefore we shall consider now  the
so called $\psi-exp$ functions [20, 7].
\end{enumerate}
} {\section {On q-extension and $\psi-$extension of higher order
$\alpha-$hyperbolic functions} In this section we shall try to
stay close to \textit{exp}. At first we shall refer to what is
known since a long time; see [12] from 1910 year and [10] for
Heine and Gauss contribution and [13] for may be application to
quantum processes description and overall theory of the so called
non-commutative geometry. Therefore we shall  consider here a
\textit{specific} example of such series $L$ which are extensions
of \textit{exp}  with  some properties surviving or being
mimicked. These are $exp_{q}$ and $exp_{\psi}$ functions. Before
doing that some:\\
 \textbf{Preliminaries:}\\
    We  perform after Heine  and Gauss [10] a replacement $x\rightarrowtail x_{q}$
    thus arriving at the standard by now deformation of the variable $x\in R$ [10,13] according to the prescription:
$$x\rightarrowtail x_{q}\equiv
\frac{1-q^{x}}{1-q}\xrightarrow[q\to 1]{} x$$

 Then  consequently we have for $n_{q}$, $q-$factorial and $q-$binomial
coefficients  ${n\choose k}_{q}\equiv\frac{n_{q}^{\underline
k}}{k_{q}!}$ where
$n_{\psi}^{\underline{k}}=n_{\psi}(n-1)_{\psi}(n-2)_{\psi}\dots(n-k+1)_{\psi}$.

Also integration and derivation [12] might be q-extended. Here we
introduce only - what is called - Jackson`s derivative
$\partial_{q}-$ a kind of difference operator.
\begin{defn} \textit{Jackson`s derivative $\partial_{q}$ is defined as
follows. Let $\phi$ denote any Laurent series. Then
$(\partial_{q}\phi)(x)=\frac{\phi(x)-\phi(qx)}{(1-q)x}$.}
\end{defn}
Naturally
 $\partial_{q}\rightarrow[q\to 1] \frac{d}{dx}$
 and is a mere of exercise to prove that $Q-$Leibniz rule holds.

\textbf{Observation 4.1.} Let $f$, $g$, $\phi$ denote Laurent
series. Let $(Q\phi)(z):=\phi(qz)$. Then $\partial_{q}(f\cdot
g)(\partial_{q}f)\cdot g+(Qf)\cdot(\partial_{q}g)$.\\ It is a easy
to see $\partial_{q}x^{n}=n_{q}x^{n-1}$ and
$\partial_{q}\exp_{q}=\exp_{q}$; $\exp_{q}[z]_{z=0}=1$ where
$q-exp$ function is defined by
$\exp_{q}[z]:=\sum_{k=0}^{\infty}\frac{z^{n}}{n_{q}!}$. Applying
now projection operators $\{\Pi_{l}\}_{l\in Z_{n}}$ to $exp_{q}$
function we get the family $\{h_{q,s}(z)\}_{s\in Z_{n}}$ of
$q-$extended hyperbolic functions of order $n$.
\begin{defn} $\{h_{q,s}(z)\}_{s\in Z_{n}}$ are
defined by
\begin{equation}\label{e24}
  h_{q,s}=\Pi_{s}\exp_{q};\quad s\in Z_{n};\quad h_{q,s}
  h_{s};\quad s\in Z_{n}.
\end{equation}
\end{defn}

Many formulas and identities $q-$extend almost automatically from
the $q=1$ case as for example those from [6] with $q-$extended
Laguerre polynomials $L_{n,q}^{(\alpha=-1)}(x)\equiv L_{n,q}(x)$
replacing the standard ones which are the so called binomial type
or convolution type depending on $n-$dependent factor. These
$q-$identities yield automatically the corresponding ones for
projected out functions $L_{l}^{\alpha}:=\Pi_{l}^{(\alpha)}L$; for
example for $L=L_{n,q}(x)$.
\\
\textbf{Example  4.1.}  (see [20]) As an example of $q-$extended
polynomial sequences we present now  the $q-$extended Laguerre
polynomials $L_{n,q}^{(\alpha=-1)}(x)\equiv L_{n,q}(x)$.
$L_{n,q}(x)=\frac{n_{q}}{n}\sum
\limits_{k=1}\limits^{n}(-1)^{k}\frac{n_{q}!}{k_{q}!}{{n-1}\choose
{k-1}}_{q}\frac{k}{k_{q}}x^{k}$ form the so called  \textit{basic
polynomial sequence} $\{L_{n,q}(X)\}_{n\geq0}$ of the operator
$Q(\partial_{q})=-\sum\limits_{k=0}\limits^{\infty}\partial_{q}^{k+1}\equiv
\frac{\partial_{q}}{\partial_{q}-1}\equiv
-[\partial_{q}+\partial_{q}^{2}+\partial_{q}^{3}+\dots]$ which is
equivalent to say that for polynomial sequence
$p_{n}(x)=L_{n,q}(x)$; $\deg p_{n}(x)=n$ the following
requirements are fulfilled: $p_{0}(x)=1$; b) $p_{n}(0)=0$; and c)
$Q(\partial_{q})p_{n}=n_{p}p_{n-1}$. One may show that the so
called q-binomiality identity holds [20]:
$$p_{n}(x+_{q}y)=\sum\limits_{k\geq 0}{n\choose
  k}_{q}p_{k}(x)p_{n-k}(x)$$
 where $E^{a}(\partial_{q})=\exp_{q}\{a\partial_{q}\}=\sum\limits_{k=0}
 \limits^{\infty}\frac {a^{k}}{k_{q}!}\partial_{q}^{k}$
  and $E^{y}(\partial_{q})p_{n}(x)\equiv p_{n}(x+_{q}y)$
The above $q-$binomial identity yield automatically the
corresponding ones for projected out new  special $q-$polynomials
$L_{l}^{\alpha}:=\Pi_{l}^{(\alpha)}L$; $L=L_{n,q}(x)$.\\ The same
could be applied also to q-extensions of  the well known
hyperbolic functions of any order. Before considering this in the
next example let us at first note that similarly to
$\frac{dh_{k}(x)}{dx}=h_{k-1}(x)$; $k\in \textbf{Z}_{n}$ from
which it follows that
$\frac{d^{k}h_{l}(x)}{dx^{k}}=h_{l\dot{-}k}(x)$; $k,l\in
\textbf{Z}_{n}$ - also the following holds.
\\
\textbf{Observation 4.2.}
\begin{equation}\label{e25}
  \partial_{q}^{k}h_{q,l}^{\alpha}=\prod\limits_{s=0}\limits^{k\dot{-}1}(1+(\alpha-1)\delta_{0,l\dot{-}s})
  h_{q,l\dot{-}k}^{\alpha}; \ k,l \in \textbf{Z}_{n}
\end{equation}

where
\begin{equation}\label{e26}
  h_{q,s}^{\alpha}\equiv \Pi_{s}^{\alpha}\exp_{q}; \hspace{1mm}
  s\in \textbf{Z}_{n}
\end{equation}

\textbf{Example 4.2.} (versus  Example  3.1.)\\
Let $\{h_{q,s}^{\alpha}(z)\}_{s\in \textbf{Z}_{n}}$ where
$h_{q,s}^{\alpha}=\Pi_{s}^{\alpha}\exp_{q}$ then
$h_{q,s}^{\alpha}(z)=\sum\limits_{k\geq0}
\frac{\alpha^{k}z^{nk+s}}{(nk+s)_{q}!}$ and $\Omega
h_{q,s}^{\alpha}=\omega^{s}h_{q,s}^{\alpha}$; $s\in
\textbf{Z}_{n}$. We shall call: $h_{q,l}^{\alpha}$ the
$l-\alpha-q-$hyperbolic series. Of course

\begin{equation}\label{e27}
  h_{q,s}^{\alpha}(z)=\frac{1}{n}\alpha^{\frac{-s}{n}}\sum\limits_{k\in
  \textbf{Z}_{n}}\omega^{-ks}\exp_{q}(\omega^{k}\sqrt[n]{\alpha}z);\hspace{1mm}s\in
  \textbf{Z}_{n}
\end{equation}
Note also that $(h_{q,j}\equiv h_{q,l}^{\alpha=1})$;
$\exp_{q}(\omega^{l}z)=\frac{1}{n}\sum\limits_{k\in
\textbf{Z}_{n}}\omega^{kl}h_{q,k}(z)$; $l\in \textbf{Z}_{n}$. The
"$\omega-$with" rescaling operator $\Omega$ has much more of
eigenvectors apart from the family represented by
\begin{equation}\label{e28}
  \Omega h_{q,s}^{\alpha}=\omega^{s}h_{q,s}^{\alpha};\hspace{1mm}s\in
  \textbf{Z}_{n}\ or\ by\ \Omega L_{l}^{\alpha}=\omega^{l}L_{l}^{\alpha}
\end{equation}
where $L_{l}^{\alpha}$ are the $l-\alpha-$Laurent series (see:
Example 3.1.). Namely: consider the generalized factorial
$n_{\psi}!\equiv n_{\psi}(n-1)_{\psi}(n-2)_{\psi}\dotsm
2_{\psi}1_{\psi}$; $0_{\psi}!=1$ for an arbitrary sequence
$\psi=\{\psi_{n}\}_{n\geq 1}$ with the condition, $\psi_{n}\neq
0$, $n\in N$. Here $n_{\psi}$ denotes the $\psi-$deformed number
where in conformity with Viskov [25] notation $n_{\psi}\equiv
\psi_{n-1}(q)\psi_{n}^{-1}(q)$ or equivalently $n_{\psi}!\equiv
\psi_{n}^{-1}(q)$ [20]. One may now define linear operator
$\partial_{\psi}$ named $\psi-$derivative on - say - polynomials
according to: $\partial_{\psi}x^{n}=n_{\psi}x^{n-1}$; $n>0$,
$\partial_{\psi}const=0$. One defines then $\psi-exp$ function
$\exp_{\psi}[z]:=\sum\limits_{k=0}\limits^{\infty}\frac{z^{n}}{n_{\psi}!}$
so that all other constructions and statements of this section
"$\psi-$extend" automatically. An so
\begin{equation}\label{e29}
  \Omega h_{\psi,s}^{\alpha}=\omega^{s}h_{\psi,s}^{\alpha};\hspace{1mm}
  s\in \textbf{Z}_{n}
\end{equation}

with self-explanatory notation:
$h_{\psi,s}^{\alpha}=\Pi_{s}^{\alpha}\exp_{\psi}$.\\
\textbf{Remark 4.1.} (see: [25] and  [20])\\ We may introduce now
${n\choose k}_{\psi}\equiv \frac{n_{\psi}^{k}}{k_{\psi}!}$ where
$n_{\psi}^{k}=n_{\psi}(n-1)_{\psi}\dots(n-k+1)_{\psi}$ and extend
a very important notion of the polynomial sequence of binomial
type. Here are examples: take for polynomial sequence
$\{p_{n}\}_{0}^{\infty}$; $\deg p_{n}=n$; $p_{n}(x)=x^{n}$ or take
$p_{n}(x)=x^{n}=x(x-1)\dots(x-n+1)$. Then one easily checks that
the following identity holds:

\begin{equation}\label{e30}
  p_{n}(x+y)\equiv \sum\limits_{k\geq 0}{n\choose
  k}_{\psi}p_{k}(x)p_{n-k}(y)
\end{equation}

Polynomial sequences satisfying  (30) are polynomial sequences of
binomial type. Polynomial sequence  $\{p_{n}\}_{0}^{\infty}$ is
then of $\psi-$binomial type if it satisfies the recurrence
$$E^{y}(\partial_{\psi})p_{n}(x)\equiv \sum\limits_{k\geq 0}{n\choose
  k}_{\psi}p_{k}(x)p_{n-k}(y)$$
where  $E^{y}(\partial_{\psi})\equiv
\exp_{\psi}\{y\partial_{\psi}\}=\sum\limits_{k=0}\limits^{\infty}
\frac{y^{k}\partial_{\psi}^{k}}{k_{\psi}!}$ is a generalized
translation operator [20]. Polynomials encompassing those of
$\psi-$binomial type are the so called [20] Sheffer
$\psi-$polynomials. In [25] (1975) - Proposition 8 - Viskov have
proved that polynomial sequence $\{p_{n}\}_{0}^{\infty}$ is
Sheffer $\psi-$polynomial if and only if its "$\psi-$generating
function" is of the form:

\begin{equation}\label{e31}
  \sum\limits_{n\geq 0}\psi_{n}p_{n}(x)z^{n}=A(z)\psi(xg(z));
\end{equation}
\begin{equation}\label{e32}
  \psi(z)=\sum\limits_{n\geq 0}\psi_{n}z^{n};\ \psi_{n}\neq 0;\
  n=0,1,2,\dots
\end{equation}

where $A(z)$, $g(z)/z$ are formal series with constant terms
different from zero. In the very important reference [7] Y.Ben
Cheikh has given important examples of decomposition of the
Boas-Buck polynomials with respect to the cyclic group $Z_{n}$. In
our notation [20] adapted to [25]  these are Sheffer
$\psi-$polynomials including polynomial sequences  of
$\psi-$binomial type.\\ \textbf{Example 4.3.} It is easy to check
that for $\psi_{n}(q)=[R(q^{n})!]^{-1}$ and $R(x)=\frac{1-x}{1-q}$
we get $\psi_{n}(q)=n_{q}$. In [25] (1975 Proposition 4) Viskov
have proved also that polynomial sequence $\{p_{n}\}_{0}^{\infty}$
is of $\psi-$binomial type if and only if its "$\psi-$generating
function" is of the form
\begin{equation}\label{e33}
  \sum\limits_{n\geq 0}\psi_{n}p_{n}(x)z^{n}=\exp_{\psi}(xg(z))
\end{equation}
for formal series  $g$ inverse to appropriate formal series (see:
[25] (1975)).  Now   for $\psi_{n}(q)=[n_{q}!]^{-1}$,
$\psi(z)=\exp_{q}\{z\}$ and "$\exp_{q}$ generating function"
takes the form
\begin{equation}\label{e34}
  \sum\limits_{n\geq
  0}\frac{z^{n}}{n_{q}!}p_{n}(x)=\exp_{q}(xg(z))
\end{equation}

If one denotes by $p_{n,s}^{\alpha,\psi}$ the following
eigenpolynomials of $\Omega$:
$p_{n,s}^{\alpha,\psi}=\Pi_{s}^{\alpha}p_{n}$; $s\in
\textbf{Z}_{n}$ and if $A(z)=1$ then for generating functions of
these special polynomials we get from [25] the following
expressions:
\begin{equation}\label{e35}
  \sum\limits_{n\geq
  0}\psi_{n}p_{n,s}^{\alpha,\psi}(x)z^{n}=h_{\psi,s}^{\alpha}(xg(z)),\
  s\in \textbf{Z}_{r}
\end{equation}
If in addition $g=id$ then
 \begin{equation}\label{e36}
  \sum\limits_{n\geq
  0}\psi_{n}p_{n,s}^{\alpha,\psi}(x)z^{n}=h_{\psi,s}^{\alpha}(xz),\
  s\in \textbf{Z}_{r}
\end{equation}
We call functions $h_{\psi,s}^{\alpha}$ the $\psi-$hyperbolic
functions. Naturally $\Omega$ $\omega-$rescales  $x$ argument of
$p_{n,s}^{\alpha,\psi}$ and $h_{\psi,s}^{\alpha}(xg(z))$; $s\in
\textbf{Z}_{r}$ and both immense sets of these special
 functions are $\omega^{s}-$homogeneous, $s\in \textbf{Z}_{r}$, which
 is equivalent to say that these are eigenfunctions of  scaling operator $\Omega$
corresponding to the eigenvalue $\omega^{s}$; $s\in
\textbf{Z}_{r}$. For $\psi_{n}(q)=n_{q}$ one gets from {e34}
$q-$deformed $\omega^{s}-$homogeneous special $q-$hyperbolic
functions $h_{q,s}^{\alpha}$ and special $\omega^{s}-$homogeneous
$q-$deformed polynomials $p_{n,s}^{\alpha,q}$. In the limit case
of $q=1$ we end up with classical special polynomials - for
example with Laguerre polynomials  [6] - and other polynomial
sequences - for example of binomial type. \\\textbf{Remark 4.2.}\\
Note that in the case of analytic functions instead of $f(x)$ one
may consider also functions with matrix arguments $f(A)$; $A\in
M_{kxk}(C)$ or arguments from associative algebras with unity over
$C$ equipped with norm in order to assure the possibility of
convergence. Hyperbolic mappings of such type might be now equally
well investigated. \\\textbf{Acknowledgements}\\
    The authors are very much indebted to Referees whose indications
allowed preparing the paper in a hopefully more desirable form.

\end{document}